\documentclass[journal]{IEEEtran}
\usepackage{cite}
\usepackage{times,amsmath,amssymb,psfrag}
\usepackage{graphicx}
\usepackage{url}
\usepackage{bm}
\usepackage{multicol}
\graphicspath{{./pics/}} \DeclareGraphicsExtensions{.pdf,.jpg,.png}
\usepackage{color}
\usepackage{hyperref}
\usepackage{epsfig}
\usepackage{float}
\usepackage{array}
\usepackage{multirow}
\usepackage{etoolbox}
\usepackage{nomencl}
\usepackage{algorithm}
\usepackage{algorithmicx}
\usepackage{algpseudocode}
\usepackage{amsmath}
\newcommand{\ct}{\textcolor{red}}

\newcommand\Nom[3][X]{\nomenclature[#1#3]{#2}{#3}}
\renewcommand\nomgroup[1]{%
  \item[\normalsize\itshape\bfseries
  \ifstrequal{#1}{I}{Indices and Sets}{%
  \ifstrequal{#1}{P}{Parameters}{%
  \ifstrequal{#1}{V}{Variables}{%
  \ifstrequal{#1}{X}{Other Symbols}{}}}}]%
}

\makenomenclature
\begin{document}
\title{Transmission System Resilience Enhancement with Extended Steady-state Security Region in Consideration of Uncertain Topology Changes}
\author{
Chong Wang,~\IEEEmembership{Member,~IEEE}, Feng Wu,~\IEEEmembership{Member,~IEEE}, Ping Ju,~\IEEEmembership{Senior Member,~IEEE},  \\ Shunbo Lei,~\IEEEmembership{Member,~IEEE}, Tianguang Lu,~\IEEEmembership{Member,~IEEE},  Yunhe Hou,~\IEEEmembership{Senior Member,~IEEE}


\thanks{
	

	
C. Wang, F. Wu and P. Ju are with the College of Energy and Electrical Engineering, Hohai University, Nanjing 211100, China (e-mail: chongwang@hhu.edu.cn, wufeng@hhu.edu.cn, pju@hhu.edu.cn).

S. Lei is with the Department of Electrical Engineering and Computer Science, University of Michigan, Ann Arbor, MI, 48109 USA (email: shunbol@umich.edu).

T. Lu is with the School of Engineering and Applied Sciences and Harvard China Project, Harvard University, Cambridge, MA 02138, USA (email: tlu@seas.harvard.edu)

Y. Hou is with the Department of Electrical and Electronic Engineering, The University of Hong Kong, Hong Kong (e-mail: yhhou@eee.hku.hk).



}
}

\date{}
\maketitle
\begin{abstract}
The increasing extreme weather events poses unprecedented challenges on power system operation because of their uncertain and sequential impacts on power systems. This paper proposes the concept of an extended steady-state security region (ESSR), and resilience enhancement for transmission systems based on ESSR in consideration of uncertain varying topology changes caused by the extreme weather events is implemented. ESSR is a ploytope describing a region, in which the operating points are within the operating constraints. In consideration of uncertain varying topology changes with ESSR, the resilience enhancement problem is built as a bilevel programming optimization model, in which the system operators deploy the optimal strategy against the most threatening scenario caused by the extreme weather events. To avoid the curse of dimensionality with regard to system topologies for a large-scale system, the Monte Carlo method is used to generate uncertain system topologies, and a recursive McCormick envelope-based approach is proposed to connect generated system topologies to optimization variables. Karush–Kuhn–Tucker (KKT) conditions are used to transform the sub-optimization model in the second level into a group of equivalent constraints in the first level. A simple test system and IEEE 118-bus system are used to validate the proposed.

\end{abstract}
\IEEEpeerreviewmaketitle
\setlength{\nomitemsep}{0.18cm}
\printnomenclature[2.0cm]

\Nom[I01]{$b,b'$}{Index of terminal buses of a line.}
\Nom[I02]{$g$}{Index of generators.}
\Nom[I03]{$k$}{Index of sequential system topology scenarios.}
\Nom[I04]{$l$}{Index of lines.}
\Nom[I05]{$t$}{Index of time periods.}
\Nom[I06]{$\mathcal{G}_b$}{Set of units at power bus $b$.} 
\Nom[I07]{$\mathcal{B}_b$}{Set of power buses connected to power bus $b$.}
\Nom[I08]{$\mathcal{L}$}{Set of lines that are in failure states for all sequential system topology scenarios.}
\Nom[I09]{$\mathcal{\bar L}$}{Set of lines that are not in failure states for all sequential system topology scenarios.}
\Nom[I10]{$\Pi_k$}{Set of system topologies from the first decision time period to the last decision time period for the $k^{th}$ sequential system topology scenario.}
\Nom[I11]{$\Omega$}{Set of all possible ESSRs corresponding to all sequential system topology scenarios.}
\Nom[I11]{$\Psi _1, \Psi _2, \Psi _K$}{Different extended steady-state security region.}

\Nom[P01]{${\bf{a}}, {\bf{A}}, {\bf{b}}, {\bf{B}}, {\bf{C}} $}{Coefficient matrices based on operational constraints.}
\Nom[P02]{$B_{bb'}$}{Electrical susceptance.}
\Nom[P03]{$M$}{A large number.}
\Nom[P04]{${\overline P}_{bb'} $}{Upper active power capacities of line $b-b'$.}
\Nom[P05]{${\overline P}_{bt} $}{Load at bus $b$ in the time period $t$.}
\Nom[P06]{${\underline P}_g, {\overline P}_g $}{Lower and upper active power limits of generator $g$.}
\Nom[P07]{$\underline{R}_{g}, \overline{R}_{g} $}{Ramp-up and ramp-down limits of generator $u$.}
\Nom[P08]{$\underline{\theta}_{b}, \overline{\theta}_{b} $}{Limits of voltage angle at bus $b$.}

\Nom[V01]{$P_{gt}, P_{g(t+1)}$}{Active power of generator $g$ at $t$ and $t+1$, respectively.}
\Nom[V02]{$P_{bb'tk}$}{Active power through line $b-b'$ in the time period $t$ under the sequential system topology scenario $k$.}
\Nom[V03]{$\theta_{btk}$}{Phase angle of bus $b$ in the time period $t$ under the sequential system topology scenario $k$.}
\Nom[V04]{$u_{bb'tk}$}{Binary variable to indicate the state of line $b-b'$ in the time period $t$ under the sequential system topology scenario $k$. `1' and `0' denote on-state and off-state, respectively.}

\Nom[V05]{${\bf{U}}$}{Binary variable vector representing different system topology.}
\Nom[V06]{${\bf{Y}}$}{Continuous variable vector representing system operating conditions.}
\Nom[V07]{${\bf{Y}^*}$}{Optimal value of the suboptimization model.}
\Nom[V08]{${{\bf{s}}^ + }$,${{\bf{s}}^ - }$}{Vectors of positive slack variables.}
\Nom[V09]{$({{\bf{s}}^ + })^ *$, $({{\bf{s}}^ - })^ *$}{Optimal values of the suboptimization model.}
\Nom[V10]{$L$}{Lagrangian funtion.}
\Nom[V11]{$f,F$}{Optimization objective functions.}
\Nom[V12]{$\alpha, \beta, \gamma$}{Lagrange multiplier vectors.}
\Nom[V13]{$\hat \alpha, \hat \beta, \hat \gamma$}{Lagrange multiplier diagonal matrices corresponding to $\alpha, \beta, \gamma$, respectively.}

\section{Introduction}
\IEEEPARstart{E}{xtreme} weather events with higher frequency and heavy intensity have great impacts on power system operation \cite{intensityR1, intensityR2}. The countermeasures against possible extreme weather events are usually not included in the conventional operational strategies. Many organizations, e.g., the United States National Research Council (NRC) \cite{NRC_R1}, the House of Lords in the United Kingdom\cite{HouseLord_R1}, the North American Electric Reliability Corporation (NERC) \cite{NERC_R1, NERC_R2}, and the United States Electric Power Research Institute (EPRI) \cite{EPRI_R1}, have already emphasized the importance of resilient operational strategies against extreme weather events. 

According to the framework of power system resilience, the strategies should be prepared in three stages of an extreme weather event, i.e., prior to the event, during the event, and after the event \cite{NRC_R1, NERC_R2}. \textit{Prior to the event}, it is necessary to perform preventive strategies and assessments to enhance the system capacity of stay standing or keeping operating in the face of disasters. To assess the possible outages caused by hurricanes, the regression models based on historical data are established \cite{EsM1, EsM2, EsM3}. In addition, the impacts of wildfires and extreme floods on power systems are analyzed \cite{Wildfires1, Floods1}. With the assessments, some preventive strategies, e.g., pre-hurricane restoration planning, network hardening, mobile energy storage pre-allocation, and microgrid construction, are investigated to increase the system capacity of stay standing in the face of extreme weather events. For example, \cite{Pr1} proposed a proactive resource allocation model in consideration of  restoration and repair for potential failures of devices located on the path of an upcoming hurricane. \cite{Pr2} presented conservative schedules for microgrids in consideration of the minimum number of vulnerable lines in service on the path of an approaching hurricane. \cite{Pr3} proposed a tri-level optimal hardening strategy to enhance the resilience of power distribution networks to protect against extreme weather events. \cite{Pr4} presented a model for resilient routing and scheduling of mobile power sources, which are pre-positioned to enable rapid pre-restoration for high survivability of critical loads. \textit{During the event}, it is necessary to perform real-time strategies to enhance the system capacity of managing the impacts as it unfolds. \cite{Dur1} presented an operational enhancement approach, which includes the assessment of the impact of severe weather events on power systems and a risk-based defensive islanding algorithm. The developed islanding algorithm can avoid the potential cascading failures during extreme weather events. \cite{WZY} developed a microgrid sectionalization method, which divides an interconnected system into several microgrids, to improve distribution system resilience. \cite{Dur2} proposed an integrated resilience response framework that links the preventive strategies and the real-time strategies, which include generator dispatch, topology switching and load shedding. \cite{Dur3} proposed a resilience-constrained unit commitment (RCUC) model in consideration of system operational constraints, heterogeneity of power flow distribution, and lines
forced outages. \cite{Dur4} proposed a Markov-based operation strategy to enhance system resilience during an unfolding extreme event. The possible sequential system topologies are modeled as Markov states, by which a recursive value function is developed with operation constraints and intertemporal
constraints. \textit{After the event}, it is necessary to perform restoration strategies to recover the system to a normal operating condition as quickly as possible. Microgrid construction is one of effective methods to restore loads quickly after disasters. \cite{After1} proposed a distribution system operational approach by constructing multiple microgrids with distributed generators to restore critical loads from the power outage. In consideration of the scarcity of power generation resources, \cite{After2} employed the concept of continuous operating time (COT) to determine the availability of microgrids to evaluate the service time and in consequence to restore critical load by using a chance-constrained model. To avoid subsequent outages, \cite{After3} presented the formation of adaptive multi-microgrids with mobile emergency resources as part of the critical service restoration strategy. After disasters, appropriately dispatching repair crews and mobile power sources can accelerate the load restoration. \cite{Lei11} proposed a resilient scheme for disaster recovery logistics to co-optimize distribution system restoration in consideration of dispatching repair crews and mobile power sources by using a mixed integer linear programming.

Based on the literature review, most of research studies focused on system hardening, mobile emergency resource allocation/dispatch, microgrid construction, and repair crew dispatch to enhance power system resilience for distribution systems. For transmission systems, the effects of these strategies are limited due to different characteristics between distribution systems and transmission systems. In practice, a proper operating point when an extreme weather event impacts the transmission systems is critical in consideration of operational constraints, e.g., ramping rates of generators, output limits of generators, and system transmission capacity. An inappropriate operating point may need to perform load shedding in face of uncertain sequential system topologies on the trajectory of an weather-related event in consideration of the limits of operational constraints. Therefore, this paper focuses on the investigation of the operating point for the transmission system to enhance the system resilience against the extreme weather event. The contributions of this paper are listed as follows: 1) The concept of an extended steady-state security region is proposed; 2) Base on the extended steady-state security region, resilience enhancement in consideration of uncertain  varying topology changes caused by the extreme weather events is implemented; 3) Resilience enhancement is constructed as a bilevel programming optimization model, with which an optimal operating point against the threatening scenarios caused by the sequential weather event.    

The remainder of this paper is organized as follows. Section II shows the extended steady-state security region during an unfolding event. Section III presents the bilevel optimization model of resilience enhancement and the solution method. The case studies are demonstrated in Section IV, and the work is concluded in Section V. 

\section{Extended Steady-state Security Region during an Unfolding Event}
\subsection{Extended Steady-state Security Region}
Due to the sequential time periods of an unfolding weather event, different transmission lines might be in failure in different time periods, and in consequence result in different topologies in different time periods. For example, the set of failure components in the time periods $t_1$ and $t_2$ may be \{$b_{12}$\} and \{$b_{12}$, $b_{14}$\}, respectively. This indicates that the line $b_{12}$ is out of service in $t_1$ and $t_2$, and the line $b_{14}$ is out of service in $t_2$. Due to the uncertain impacts of the unfolding event on the system, the possible failure components may have different scenarios, e.g., \{$b_{12}$\} and \{$b_{12}$, $b_{16}$\} in $t_1$ and $t_2$, respectively. For all of these possible scenarios, the system operators expect feasible operating points in each decision period. However, the operating points over sequential time periods are determined by operational constraints, e.g., generators' ramping rates, thermal capacity limits of lines, and generators' capacity limits. An appropriate operating point in the time period $t_0$ play an important role in the subsequent operating points. We use the steady-state security region (SSR) \cite{SR_FFWU} to illustrate the importance of the operating point.   
\begin{figure}[!h]
	\centering
	\includegraphics[width=7cm]{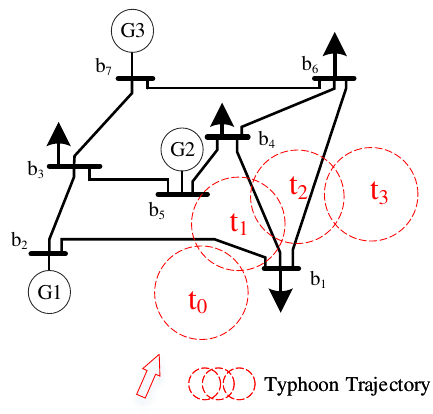}
	\caption{An example of components' failure scenarios during an unfolding event}
	\label{StateDemonstrate_fig}
\end{figure}  

\begin{figure}[!h]
	\centering
	\includegraphics[width=7cm]{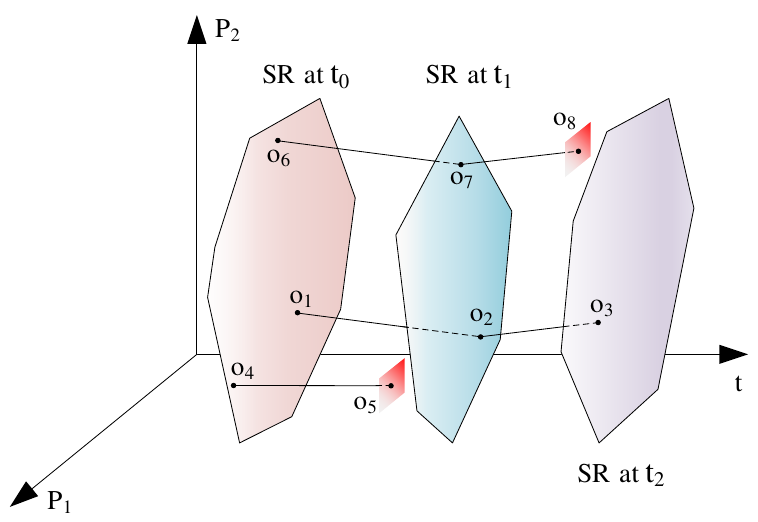}
	\caption{An example of security region (SR) on a sequential timeline}
	\label{ESR_fig}
\end{figure}

The SSR is defined as a multi-dimensional space with the set of injected power as axes, where the operational constraints work as boundaries. Based on SSR, it is easy to check whether an operating condition is steady-state secure. Take the scenario in Fig. \ref{ESR_fig} as example, in which there are two generators. The operating conditions $o_1$, $o_4$, and $o_6$ at time $t_0$ are steady-state secure. Usually, SSR is constructed to deal with one time period. However, in the practical system, strategies are usually needed over sequential time periods. When considering the impacts of sequential time periods, the scenarios will be different. For example, $o_1$ can reach $o_2$, and then can reach $o_3$ in consideration of the generators' ramp-rates in $t_1$, $t_2$, and $t_3$, respectively. However, $o_4$ can only reach the red rectangular region in $t_2$, where the corresponding operating region is not within SSR in $t_2$. Similarly, $o_7$ can only reach the red rectangular region in $t_3$, where the corresponding operating region is not within SSR in $t_3$. 
From the above cases, we conclude that the impacts of the sequential strategies in different time periods should be included. Therefore, an extended steady-state security region (ESSR) is proposed based on SSR.  

SSR is assumed to be expressed as follows.
\begin{align}
& \psi {\rm{ = \{ }}{\bf{y}}{\rm{|}}{\bf{a} \cdot \bf{y}} \le {\bf{b}}{\rm{\} }} \label{SSR}
\end{align}
where $\bf{y}$ is the variable vector, in which generators' outputs are included, and $\bf{a}$ and $\bf{b}$ are coefficient matrices representing operating constraints. Based on (\ref{SSR}), ESSR can be expressed as follows.  
\begin{subequations}
	\label{ESSR}
	\begin{align}
	& \Psi {\rm{ = \{ }}{\bf{Y}}{\rm{|}}{\bf{A \cdot Y}} \le {\bf{B}}{\rm{\} }} \label{ESSR1}\\[8pt]
	& {\bf{Y}} = [{\bf{y}}_1^T, \cdots ,{\bf{y}}_t^T \cdots ,{\bf{y}}_N^T]^T \label{ESSR2}\\[8pt]
	& {\bf{B}} = [{\bf{b}}_1^T, \cdots ,{\bf{b}}_t^T \cdots ,{\bf{b}}_N^T]^T \label{ESSR3}\\[8pt]
	& {\bf{A}} = \left[ {\begin{array}{*{20}{c}}
		{{{\bf{a}}_1}}&{{{\bf{a}}_{12}}}&{}&{}&{}\\
		{{{\bf{a}}_{21}}}&{{{\bf{a}}_2}}&{}&{}&{}\\
		{}&{}& \ddots &{}&{}\\
		{}&{}&{}&{{{\bf{a}}_{N - 1}}}&{{{\bf{a}}_{N - 1,N}}}\\
		{}&{}&{}&{{{\bf{a}}_{N,N - 1}}}&{{{\bf{a}}_N}}
		\end{array}} \right]
	\end{align}
\end{subequations}
where $\bf{Y}$ is the variable vector, in which generators' outputs in different time periods are included, and $\bf{A}$ and $\bf{B}$ are coefficient matrices representing operating constraints in different time periods. ${\bf{a}}_{12}$ and ${\bf{a}}_{21}$ represent the coupling relations between the first time period and the second time period, and they correspond to the generators' ramp-rates.
   
\subsection{ESSR with Uncertain Varying Topology Changes}
\begin{figure}[!h]
	\centering
	\includegraphics[width=9cm]{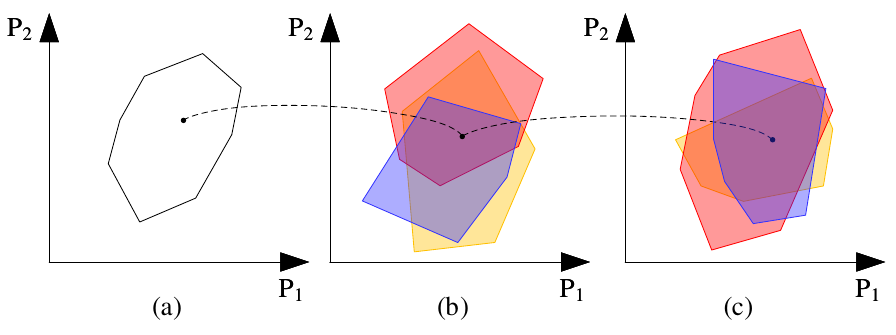}
	\caption{Extended Security region (ESR) illustrated by three 2-dimension figures}
	\label{ESR_Uncertainty_fig}
\end{figure}

ESSR addresses the issue of sequential characteristics of strategies on the system, i.e, the sequential impacts of extreme weather events on the system. However, the uncertain impacts of extreme weather events on the system should be also included. Take the scenario in Fig. \ref{ESR_Uncertainty_fig} as an example, and we illustrate the 3-dimension figure by means of three 2-dimension figures to read them easily. Fig. \ref{ESR_Uncertainty_fig} (a), (b), and (c) represent the regions in $t_0$, $t_1$, and $t_2$, respectively. When considering the strategy in $t_0$, the system topology in the future time periods $t_1$ and $t_2$ cannot be determined because the lines may be out of service due to extreme weather events, and this results in different possible ESSRs. When dispatching the system in $t_0$, it is desired to get a strategy, which can reach the feasible region in $t_1$ and $t_2$ in consideration of different possible ESSRs. Therefore, we define a new set $\Omega$ that include all possible ESSRs, and it is expressed as follows.
\begin{align}
& \Omega  = \{ {\Psi _1},{\Psi _2}, \cdots ,{\Psi _K}\} \label{Equ_OmegaRR}
\end{align}
where $\Psi _1$, $\Psi _2$, and $\Psi _K$ are different ESSRs that correspond to different sequential system topologies $\Pi _1=\{\pi_{11}, \cdots, \pi_{1N}\}$, $\Pi _2 =\{\pi_{21}, \cdots, \pi_{2N}\}$, and $\Pi _K =\{\pi_{K1}, \cdots, \pi_{KN}\}$. $\pi_{KN}$ is a system topology scenario of the $K^{th}$ ESSR at time $N$. The sequential system topologies $\Pi _1$, $\Pi _2$, and $\Pi _K$ can be generated stochastically by using the Monte Carlo method based on failure probability due to the extreme weather events \cite{FailureModel1, FailureModel2, FailureModel3} for a large scale system.   

\section{Bilevel Optimization Model of Resilience Enhancement}
This section first presents the mathematical formulation of ESSR in consideration of with uncertain varying topology changes, and then the bilevel optimization model is developed.  
\subsection{Formulation of ESSR}
For a given sequential system topology $\Pi _k$, we can list its corresponding ESSR as follows. 
\begin{subequations}
	\label{ESSR_F}
	\begin{align}
	& \sum\limits_{g \in {G_b}} {{P_{gt}}}  - {P_{bt}} + \sum\limits_{b' \in {B_b}} {{P_{bb't}}}  = 0\quad \forall b,t \label{ESSR_F1}\\[4pt]
	& {B_{bb'}}({\theta _{bt}} - {\theta _{b't}}) = {P_{bb't}}\quad \forall (b,b') \in \Pi _K,t \label{ESSR_F2}\\[4pt]
	& {{\underline P}_{bb'}} \le {P_{bb't}} \le {{\overline P}_{bb'}}\quad \forall (b,b') \in \Pi _K,t \label{ESSR_F3}\\[4pt]
	& {{\underline R}_g} \le {P_{g(t + 1)}} - {P_{gt}} \le {{\overline R}_g}\quad \forall g,t \label{ESSR_F4}\\[4pt]
	& {{\underline P}_g} \le {P_{gt}} \le {{\overline P}_g}\quad \forall g,t \label{ESSR_F5}\\[4pt]
	& {{\underline \theta}_b} \le {\theta _{bt}} \le {{\overline \theta }_b}\quad \forall b,t \label{ESSR_F6}
	\end{align}
\end{subequations} 
where (\ref{ESSR_F1}) shows the power balance constraint, (\ref{ESSR_F2}) presents the physical relations between voltage angles and power flows through lines, (\ref{ESSR_F3}) is the power limit of lines, (\ref{ESSR_F4}) is the limit of generators' ramp-rates, (\ref{ESSR_F5}) shows the limit of generators' outputs, and (\ref{ESSR_F6}) is the voltage limit. To make topology changes easily embedded in the model, ${P_{bb't}}$ is represented as an explicit function with regard to the voltage angles $\theta_b$ and $\theta_b'$ but not the injected power. (\ref{ESSR_F}) can be expressed as a generic form as follows.
\begin{align}
& \exists {\bf{Y}}:\;{\bf{A \cdot Y}} \le {\bf{B}} \label{ESSR_F_Generic}
\end{align}
The physical meaning is that there exist a strategy satisfying the operating constraints for a given ESSR. To find a strategy, we can construct a new optimization model by introducing positive slack vectors ${\bf{s}}^{+}$ and ${\bf{s}}^{-}$.
\begin{subequations}
	\label{OPT1}
	\begin{align}
	& f({\bf{Y}}) = \mathop {\min }\limits_{{\bf{Y}},{{\bf{s}}^ + },{{\bf{s}}^ - }} \;{{\bf{1}}^T}{{\bf{s}}^ + } + {{\bf{1}}^T}{{\bf{s}}^ - } \label{OPT1_1}\\[4pt]
	& s.t. \;\; {\bf{A \cdot Y}} +  {{\bf{s}}^ + } -  {{\bf{s}}^ - } \le {\bf{B}} \label{OPT1_2}\\[4pt]
	&\quad\quad {{\bf{s}}^ + } \ge {\bf{0}},\quad {{\bf{s}}^ - } \ge {\bf{0}} \label{OPT1_3}
	\end{align}
\end{subequations} 
where $\bf{1}$ and $\bf{0}$ are vectors that consist of 1 and 0 with proper dimensions, respectively. When having $f({\bf{Y}})=0$, we can conclude $\bf{Y} \ne \emptyset $. If $f({\bf{Y}})>0$, it indicates that no strategy can satisfy the operating constraints, i.e, $\bf{Y} = \emptyset$. Therefore, we can conclude that the necessary and sufficient condition for $\bf{Y} \ne \emptyset $ is  $f({\bf{Y}})=0$.

\subsection{Formulation of ESSR with uncertain varying topology changes}
An extreme weather event may lead to different system topologies on its trajectory. When considering possible topology changes, ESSR can be represented as follows. 
\begin{subequations}
	\label{ESSR_Fpt}
	\begin{align}
	& \sum\limits_{g \in {G_b}} {{P_{gt}}}  - {P_{bt}} + \sum\limits_{b' \in {B_b}} {{P_{bb'tk}}}  = 0\quad \forall b,t,k \label{ESSR_Fpt1}\\[4pt]
	& {B_{bb'}}({\theta _{btk}} - {\theta _{b'tk}}) = {P_{bb'tk}}\quad \forall (b,b') \in \Pi _k,t,k \label{ESSR_Fpt2}\\[4pt]
	& {{\underline P}_{bb'}} \le {P_{bb'tk}} \le {{\overline P}_{bb'}}\quad \forall (b,b') \in \Pi _k,t,k \label{ESSR_Fpt3}\\[4pt]
	& {{\underline R}_g} \le {P_{g(t + 1)}} - {P_{gt}} \le {{\overline R}_g}\quad \forall g,t \label{ESSR_Fpt4}\\[4pt]
	& {{\underline P}_g} \le {P_{gt}} \le {{\overline P}_g}\quad \forall g,t \label{ESSR_Fpt5}\\[4pt]
	& {{\underline \theta}_b} \le {\theta _{btk}} \le {{\overline \theta }_b}\quad \forall b,t,k \label{ESSR_Fpt6}
	\end{align}
\end{subequations}
where (\ref{ESSR_Fpt1})-(\ref{ESSR_Fpt6}) include the possible topology changes by introducing the topology scenario notation $k$. Even though (\ref{ESSR_Fpt1})-(\ref{ESSR_Fpt6}) are directly related to the system topologies, they cannot be optimized directly because they are not explicit models with regards to optimization variables. To relate optimization variables, we reformulate the model by introducing integer variables that represent on-off states of lines and in consequence represent the system topologies.   
\begin{subequations}
	\label{ESSR_OFpt}
	\begin{align}
	& \sum\limits_{g \in {G_b}} {{P_{gt}}}  - {P_{bt}} + \sum\limits_{b' \in {B_b}} {{P_{bb'tk}}}  = 0\quad \forall b,t,k \label{ESSR_OFpt1}\\[4pt]
	& \begin{array}{l}
	{B_{bb'}}({\theta _{btk}} - {\theta _{b'tk}}) - {P_{bb'tk}} + (1 - {u_{bb'tk}})M \ge 0\\
	\quad \quad \quad\quad\quad \quad \quad \quad \quad \quad \quad \forall (b,b') \in l,l \in \mathcal{L},t
	\end{array}  \label{ESSR_OFpt2}\\[4pt]
	& \begin{array}{l}
	{B_{bb'}}({\theta _{btk}} - {\theta _{b'tk}}) - {P_{bb'tk}} - (1 - {u_{bb'tk}})M \le 0\quad \\
	\quad \quad \quad\quad\quad \quad \quad \quad \quad \quad \quad \forall (b,b') \in l,l \in \mathcal{L},t
	\end{array}   \label{ESSR_OFpt3}\\[4pt]
	& {u_{bb'tk}}{{\underline P}_{bb'}} \le {P_{bb'tk}} \le {u_{bb'tk}}{{\overline P}_{bb'}}\quad \forall (b,b') \in l,l \in \mathcal{L},t  \label{ESSR_OFpt4}\\[4pt]
	& {B_{bb'}}({\theta _{btk}} - {\theta _{b'tk}}) = {P_{bb'tk}}\quad \forall (b,b') \in l,l \in \mathcal{\overline L},t,k \label{ESSR_OFpt5}\\[4pt]
	& {{\underline P}_{bb'}} \le {P_{bb'tk}} \le {{\overline P}_{bb'}}\quad \forall (b,b') \in l,l \in \mathcal{\overline L},t,k \label{ESSR_OFpt6}\\[4pt]
	& {{\underline R}_g} \le {P_{g(t + 1)}} - {P_{gt}} \le {{\overline R}_g}\quad \forall g,t \label{ESSR_OFpt7}\\[4pt]
	& {{\underline P}_g} \le {P_{gt}} \le {{\overline P}_g}\quad \forall g,t \label{ESSR_OFpt8}\\[4pt]
	& {{\underline \theta}_b} \le {\theta _{btk}} \le {{\overline \theta }_b}\quad \forall b,t,k \label{ESSR_OFpt9}
	\end{align}
\end{subequations}
where (\ref{ESSR_OFpt2})-(\ref{ESSR_OFpt4}) represent the constraints with regards to the lines that are impacted by the extreme weather event. ${u_{bb'tk}} \; \forall (b,b') \in l,l \in \mathcal{L},t$ can cover all system topologies in all possible EESRs. In practice, a large scale system will suffer from curse of dimensionality with regard to system topologies, but some system topologies caused by the extreme weather event happen with a very small probability. To avoid curse of dimensionality, the uncertain system topologies with very small probabilities can be ignored. To this end, we can use Monte Carlo to generate system topology scenarios. To include the generated system topology scenarios as variables in the model, the additional constraints should be included. We use an example to explain this. We have three lines $l_1$, $l_2$, and $l_3$. There will be $2^3=8$ topology scenarios if we consider all scenarios, and all topology scenarios can be represented by binary variables $x_1$, $x_2$, and $x_3$, which represent the states of three lines. In practice, we may just need to consider some topology scenarios, e.g., the topology scenarios $\{x_1, x_2, x_3\} \in \{ \{1, 1, 1\}, \{0, 1, 1 \}\}$. `1' and `0' represent the on-state and off-state of the lines, respectively. To model these two topology scenarios in the model, we introduce two binary variables $z_1$ and $z_2$ to represent selection of scenarios.  $z_1=1$ denotes selection of the scenario $\{1, 1, 1\}$, and $z_2=1$ denotes selection of the scenario $\{0, 1, 1\}$. In this case, the constraints with regard to the topology scenarios can be expressed as follows.          
\begin{subequations}
	\label{Exam_case}
	\begin{align}
	& {z_1} = {x_1}{x_2}{x_3} \label{Exam_case1} \\[2pt]
	& {z_2} = (1 - {x_1}){x_2}{x_3} \label{Exam_case2} \\[2pt]
	& {z_1} + {z_2} = 1 \label{Exam_case3} \\[2pt]
	& {z_1}, {z_2}, {x_1}, {x_2}, {x_3} \in \{0, 1\} \label{Exam_case4} 
	\end{align}
\end{subequations} 
where (\ref{Exam_case1}) and (\ref{Exam_case2}) represent the relations between the scenario selection and the line states. (\ref{Exam_case3}) denotes that only one scenario can be selected. The multilinear constraint (\ref{Exam_case1}) can be represented as linear constraints by using the McCormick envelope method with new binary variables $\alpha=x_1x_2$ as follows. 
\begin{subequations}
	\label{ReExpress}
	\begin{align}
	& \alpha  \ge {x_1} + {x_2} - 1,\quad \alpha  \le {x_1},\quad \alpha  \le {x_2} \label{ReExpress1} \\[2pt]
	& {z_1} \ge \alpha  + {x_3} - 1,\quad {z_1} \le \alpha ,\quad {z_1} \le {x_3} \label{ReExpress2} 
	\end{align}
\end{subequations} 
where (\ref{ReExpress1}) and (\ref{ReExpress2}) are equivalent to (\ref{Exam_case1}) because all variables in (\ref{Exam_case1}), (\ref{ReExpress1}), and (\ref{ReExpress2}) are binary. 

Similarly, the multilinear constraints (\ref{Exam_case2}) can be represented, with new binary variables $\beta=(1-x_1)x_2$, as follows.  
\begin{subequations}
	\label{ReExpressN}
	\begin{align}
	& \beta  \ge  - {x_1} + {x_2},\quad \beta  \le 1 - {x_1},\quad \beta  \le {x_2} \label{ReExpressN1} \\[2pt]
	& {z_2} \ge \beta  + {x_3} - 1,\quad {z_2} \le \beta ,\quad {z_2} \le {x_3} \label{ReExpressN2} 
	\end{align}
\end{subequations} 

Based on the above example, we can write a generic form.
\begin{subequations}
	\label{GenericN}
	\begin{align}
	& {z_k} = \prod\limits_{ (b,b') \in l,l \in L,t} {(1 - 2{r_{bb'tk}})} (1 - {u_{bb'tk}} - {r_{bb'tk}}) \quad \forall k  \label{GenericN1} \\[2pt]
	& \quad\quad\quad \sum\limits_k {{z_k}}  = 1 \label{GenericN2} 
	\end{align}
\end{subequations} 
where ${r_{bb'tk}}$ is a known value that corresponds to the element of $\{ \{1, 1, 1\}, \{0, 1, 1 \}\}$, and ${{z_k}}$ and ${u_{bb'tk}}$ are binary variables. The multilinear functions in (\ref{GenericN1}) is represented as a group of linear constraints by means of a recursive procedure with the McCormick envelope method. In (\ref{GenericN1}), the multilinear function with all binary variables is $u_1u_2\cdots u_S$, in which the subscript $S$ is the number of binary variables for one ${{{\Psi }}_k}$.  $u_1u_2\cdots u_S$ can be expressed by introducing binary variables $\zeta_2, \; \zeta_3, \; \cdots, \zeta_{S}$ as follows.
\begin{align}
& \begin{array}{l}
{\zeta _2} = {u_1}{u_2}\\
{\zeta _3} = {\zeta _2}{u_3}\\
\cdots \\
{\zeta _S} = {\zeta _{S - 1}}{u_S}
\end{array} \label{NewCons} 
\end{align}
With the McCormick envelope method, (\ref{NewCons}) can be equivalent to the following constraints.
\begin{align}
\begin{array}{l}
{\zeta_{2}} \ge {u _{2}} + {u _{1}} - 1\\
{\zeta_{2}} \le {u _{1}}\\
{\zeta_{i}} \ge 0\quad (i = 2, \cdots ,S)\\
{\zeta_{i}} \le {u _{i}}\quad (i = 2, \cdots ,S)\\
{\zeta_{i}} \ge {u _{i}} + {\zeta_{i - 1}} - 1\quad (i = 3, \cdots ,S)\\
{\zeta_{i}} \le {\zeta_{i - 1}}\quad (i = 3, \cdots ,S)
\end{array}
\label{NewConstraints}
\end{align}

Associated with (\ref{GenericN2})-(\ref{NewConstraints}), the optimization model (\ref{ESSR_OFpt}),  can be expressed as a generic form as follows.
\begin{align}
& \forall {\bf{U}}, \;\exists {\bf{Y}}:\;{\bf{A \cdot Y}} + {\bf{C \cdot U}} \le {\bf{B}} \label{ESSR_Fpt_Generic}
\end{align}
where ${\bf{U}}$ denotes the vector of binary variables in the model. ${\bf{Y}}$ is the vector of continuous variables representing operating conditions. ${\bf{A}}$, ${\bf{B}}$, and ${\bf{C}}$ are the matrices corresponding to coefficients in the constraints (\ref{ESSR_OFpt1})-(\ref{ESSR_OFpt9}). The physical meaning of (\ref{ESSR_OFpt}) is that there is a strategy satisfying the operating constraints for all ESSRs. Based on the optimization (\ref{OPT1}), we can model the problem as follows.
\begin{subequations}
	\label{OPT2}
	\begin{align}
	& F({\bf{Y}}) = \mathop {\min }\limits_{{\bf{Y}},{{\bf{s}}^ + },{{\bf{s}}^ - }} \;{{\bf{1}}^T} \cdot {{\bf{s}}^ + } + {{\bf{1}}^T} \cdot {{\bf{s}}^ - } \label{OPT2_1}\\[2pt]
	& s.t. \;\; \forall {\bf{U}}, \;\exists {\bf{Y}}, {{\bf{s}}^ + } \ge {\bf{0}}, {{\bf{s}}^ - } \ge {\bf{0}} \label{OPT2_2}\\[2pt]
	&\quad\quad {\bf{A \cdot Y}}+ {\bf{C \cdot U}} + {{\bf{s}}^ + } - {{\bf{s}}^ - } \le {\bf{B}} \label{OPT2_3}
	\end{align}
\end{subequations}
where the optimization model (\ref{OPT2}) is to find a strategy that satisfies the worst-case scenario. Therefore, (\ref{OPT2}) can be rewritten as a linear max-min optimization model.
\begin{subequations}
	\label{OPT3}
	\begin{align}
	& F({\bf{Y}}) = \mathop {\max }\limits_{\bf{U} } \mathop {\min }\limits_{{\bf{Y}},{{\bf{s}}^ + },{{\bf{s}}^ - }} \;{{\bf{1}}^T}\cdot{{\bf{s}}^ + } + {{\bf{1}}^T}\cdot{{\bf{s}}^ - } \label{OPT3_1}\\[3pt]
	& s.t. \;\; {\bf{A \cdot Y}}+ {\bf{C \cdot U}} + {{\bf{s}}^ + } - {{\bf{s}}^ - } \le {\bf{B}} \label{OPT3_2}\\[3pt]
	&\quad\quad {{\bf{s}}^ + } \ge {\bf{0}},\quad {{\bf{s}}^ - } \ge {\bf{0}} \label{OPT3_3}
	\end{align}
\end{subequations} 
where (\ref{OPT3}) can be further reformulated as a bilevel optimization model as follows.
\begin{subequations}
	\label{OPT4}
	\begin{align}
	& F({\bf{Y}}) = \mathop {\max }\limits_{ \bf{U} } \;{{\bf{1}}^T}\cdot{({{\bf{s}}^ + })^ * } + {{\bf{1}}^T}\cdot{({{\bf{s}}^ - })^ * } \label{OPT4_1}\\[4pt]
	& s.t. \;\; \mathop {\min }\limits_{{\bf{Y}},{{\bf{s}}^ + },{{\bf{s}}^ - }} \;{{\bf{1}}^T}\cdot{{\bf{s}}^ + } + {{\bf{1}}^T}\cdot{{\bf{s}}^ - } \label{OPT4_2} \\
	&\quad\quad\quad s.t. \quad {\bf{A \cdot Y}} + {\bf{C \cdot U}}+ {{\bf{s}}^ + } - {{\bf{s}}^ - } \le {\bf{B}} \label{OPT4_3} \\
	&\quad\quad\quad\quad\quad\; {{\bf{s}}^ + } \ge {\bf{0}} \label{OPT4_4} \\
	&\quad\quad\quad\quad\quad\; {{\bf{s}}^ - } \ge {\bf{0}} \label{OPT4_5}
	\end{align}
\end{subequations} 
where $({{\bf{s}}^ + })^ *$ and $({{\bf{s}}^ - })^ *$ are the optimal solutions of the sub-optimization model (\ref{OPT4_2})-(\ref{OPT4_4}).

\subsection{Solution}
Because the sub-optimization model (\ref{OPT4_2})-(\ref{OPT4_5}) is the constraint of the main model, it is not easy to solve the model directly. The sub-optimization model (\ref{OPT4_2})-(\ref{OPT4_5}) can be equivalent to several inequality constraints and equality constraints by using Karush–Kuhn–Tucker (KKT) conditions with new variable vectors ${\bf{\alpha}}$, ${\bf{\beta}}$ and ${\bf{\gamma}}$ corresponding to (\ref{OPT4_3}), (\ref{OPT4_4}), and (\ref{OPT4_5}), respectively. The lagrangian corresponding to the sub-optimization model (\ref{OPT4_2})-(\ref{OPT4_5}) can be expressed as follows.
\begin{align}
& \begin{array}{l}
L = {{\bf{1}}^T}{{\bf{s}}^ + } + {{\bf{1}}^T}{{\bf{s}}^ - } + \\[4pt]
\quad \;\;{{\bf{\alpha }}^T}({\bf{AY}} + {\bf{CU}} + {{\bf{s}}^ + } - {{\bf{s}}^ - } - {\bf{B}}) - \\[4pt]
\quad \;\;{{\bf{\beta }}^T}{{\bf{s}}^ + } - {{\bf{\gamma }}^T}{{\bf{s}}^ - }
\end{array} \label{lagrangian_equ}
\end{align}
Based on the lagrangian, $\bf{Y^*}$, ${({{\bf{s}}^ + })^ * }$, and ${({{\bf{s}}^ - })^ * }$ of the sub-optimization model (\ref{OPT4_2})-(\ref{OPT4_5}) satisfy the KKT conditions that are be expressed as follows.
\begin{subequations}
	\label{OPT5}
	\begin{align}
	& \frac{{\partial L}}{{\partial {\bf{Y}}}} = {{\bf{A}}^T} \cdot {\bf{\alpha }} = 0 \label{OPT5_1}\\[4pt]
	& \frac{{\partial L}}{{\partial {{\bf{s}}^ + }}} = {\bf{1}} + {\bf{\alpha }} - {\bf{\beta }} = 0 \label{OPT5_2} \\
	& \frac{{\partial L}}{{\partial {{\bf{s}}^ - }}} = {\bf{1}} - {\bf{\alpha }} - {\bf{\gamma }} = 0 \label{OPT5_3} \\
	& {\bf{A}} \cdot {{\bf{Y}}^*} + {\bf{C}} \cdot {\bf{U}} + {({{\bf{s}}^ + })^ * } - {({{\bf{s}}^ - })^ * } - {\bf{B}} \le {\bf{0}} \label{OPT5_4} \\
	& - {({{\bf{s}}^ + })^ * } \le {\bf{0}} \label{OPT5_5} \\
	& - {({{\bf{s}}^ - })^ * } \le {\bf{0}} \label{OPT5_6} \\
	&  {\bf{\alpha }} \ge {\bf{0}} \label{OPT5_7} \\
	&  {\bf{\beta }} \ge {\bf{0}} \label{OPT5_8} \\
	&  {\bf{\gamma }} \ge {\bf{0}} \label{OPT5_9} \\
	& {\bf{\hat \alpha }} \cdot ({\bf{A}} \cdot {{\bf{Y}}^*} + {\bf{C}} \cdot {\bf{U}} + {({{\bf{s}}^ + })^ * } - {({{\bf{s}}^ - })^ * } - {\bf{B}}) = 0 \label{OPT5_10} \\
	& {\bf{\hat \beta }} \cdot {({{\bf{s}}^ + })^ * } = {\bf{0}} \label{OPT5_11} \\
	& {\bf{\hat \gamma }} \cdot {({{\bf{s}}^ - })^ * } = {\bf{0}} \label{OPT5_12} 
	\end{align}
\end{subequations} 
where ${\bf{\hat \alpha }}$, ${\bf{\hat \beta }}$, and ${\bf{\hat \gamma }}$ are diagonal matrices corresponding to the vectors ${\bf{\alpha }}$, ${\bf{\beta }}$, and ${\bf{\gamma }}$, respectively. 

Based on the equivalent constraints, the optimization model can be expressed as 
\begin{subequations}
	\label{OPT6}
	\begin{align}
	& \mathop {\max }\limits_{{\bf{U}},{{\bf{Y}}^*},{{({{\bf{s}}^ + })}^ * },{{({{\bf{s}}^ - })}^ * }} \;\;\;{{\bf{1}}^T} \cdot {({{\bf{s}}^ + })^ * } + {{\bf{1}}^T} \cdot {({{\bf{s}}^ - })^ * } \label{OPT6_1}\\[4pt]
	& \quad\quad\quad s.t. \quad\quad\quad\quad\quad   (\ref{OPT5_1})-(\ref{OPT5_12})  \label{OPT6_2}
	\end{align}
\end{subequations} 
where (\ref{OPT6}) is a mixed integer linear programming model, which can be solved by the solvers like CPLEX and Gurobi. If the maximum value of the model (\ref{OPT6}) is not $0$, it indicates that load shedding should be performed to ensure the existence of feasible operating points satisfying generated topology scenarios in each time period. If the maximum value of the model (\ref{OPT6}) is $0$, there exist operating points without load shedding in each time period, which satisfy generated  topology scenarios due to the extreme weather event. With the extended steady-state security region-based model, we can find an optiaml strategy to enhance the system resilience against the extreme weather event.

\section{Case Studies}
Two test systems, a simple system and a revised IEEE 118-bus system, are employed to verify the proposed model and the algorithm. The cases are tested in MATLAB 2017a using CPLEX 12.6 on a computer with 3.1 GHz i7 processors and 16 GB RAMS.  

\subsection{Case 1: A Simple System}
This simple system is used to illustrate the proposed model with detailed results. The system topology and the typhnoon trajectory are shown in Fig. \ref{StateDemonstrate_fig}. There are three generators $G_1$, $G_2$, and $G_3$, which are connected to the buses $b_2$, $b_5$, and $b_7$, respectively. Loads connected to the buses $b_1$, $b_3$, $b_4$, and $b_6$ are $0.4$, $0.4$, $0.6$, and $0.6$ (p.u.), respectively. The lower/upper limits of each generator are 0.2 (p.u.) and 2.5 (p.u.), respectively. As shown in Fig. \ref{StateDemonstrate_fig}, there are four decision epoches on the typhnoon trajectory. In $t_1$ and $t_2$, the system is directly impacted by the typhnoon. Even though the system is not directly impacted by the typhnoon in $t_0$, the dispatch strategy in $t_0$ has impacts on the following dispatch strategies in $t_1$ and $t_2$ due to the operational constraints, the generators' ramping rates, and possible topology changes caused by the typhnoon. Therefore, we focus on the dispatch strategies in $t_0$, $t_1$, and $t_2$. 
\begin{table}[!h]
	\centering
	\renewcommand{\arraystretch}{1.2}
	\caption{Line Data}
	\label{table1}
	\begin{tabular}{p{0.6cm}<{\centering}p{1.5cm}<{\centering}p{1.5cm}<{\centering}p{2.5cm}<{\centering}}
		\hline
		No. & From Bus            & To Bus  & Line Capacity (p.u.)   \\ \hline
		1   & 1                   & 2       & 0.65  \\ 
		2   & 1                   & 4       & 0.65  \\ 
		3   & 1                   & 6       & 0.65  \\ 
		4   & 2                   & 3       & 0.9  \\ 
		5   & 3                   & 5       & 0.9  \\ 
		6   & 4                   & 5       & 0.8  \\ 
		7   & 4                   & 6       & 0.8  \\ 
		8   & 6                   & 7       & 0.9  \\ 
	    9   & 3                   & 7       & 0.9 \\ \hline
	\end{tabular}
\end{table} 
\begin{table}[!h]
	\renewcommand{\arraystretch}{1.2}
	\caption{Component Failure Scenarios}
	\label{IEEE7_PossibleTopo}
	\begin{tabular}{p{2.2cm}<{\centering}p{2.2cm}<{\centering}p{2.5cm}<{\centering}}
		\hline
		\multirow{2}{*}{No.} & \multicolumn{2}{c}{Component Failure}   \\ \cline{2-3} 
		& $t_1$                                                            & $t_2$                  \\ \hline
		1                  & -                                                                & -                      \\ 
		2                  & -                                                                & $b_1$-$b_4$                       \\ 
		3                  & -                                                                & $b_1$-$b_6$                 \\ 
		4                  & -                                                                & $b_1$-$b_4$, $b_1$-$b_6$                           \\ 
		5                  & $b_1$-$b_2$                                                      & -                          \\ 
		6                  & $b_1$-$b_2$                                                      & $b_1$-$b_4$                          \\ 
		7                  & $b_1$-$b_2$                                                      & $b_1$-$b_6$                           \\ 
		8                  & $b_1$-$b_2$                                                      & $b_1$-$b_4$, $b_1$-$b_6$                          \\ 
		9                  & $b_1$-$b_4$                                                      & -                          \\ 
		10                 & $b_1$-$b_4$                                                      & $b_1$-$b_6$                          \\ 
		11                 & $b_1$-$b_2$, $b_1$-$b_4$                                         & -                         \\ 
		12                 & $b_1$-$b_2$, $b_1$-$b_4$                                         & $b_1$-$b_6$ \\  \hline
	\end{tabular}
\end{table}

\begin{figure}[!h]
	\centering
	\includegraphics[width=5.5cm]{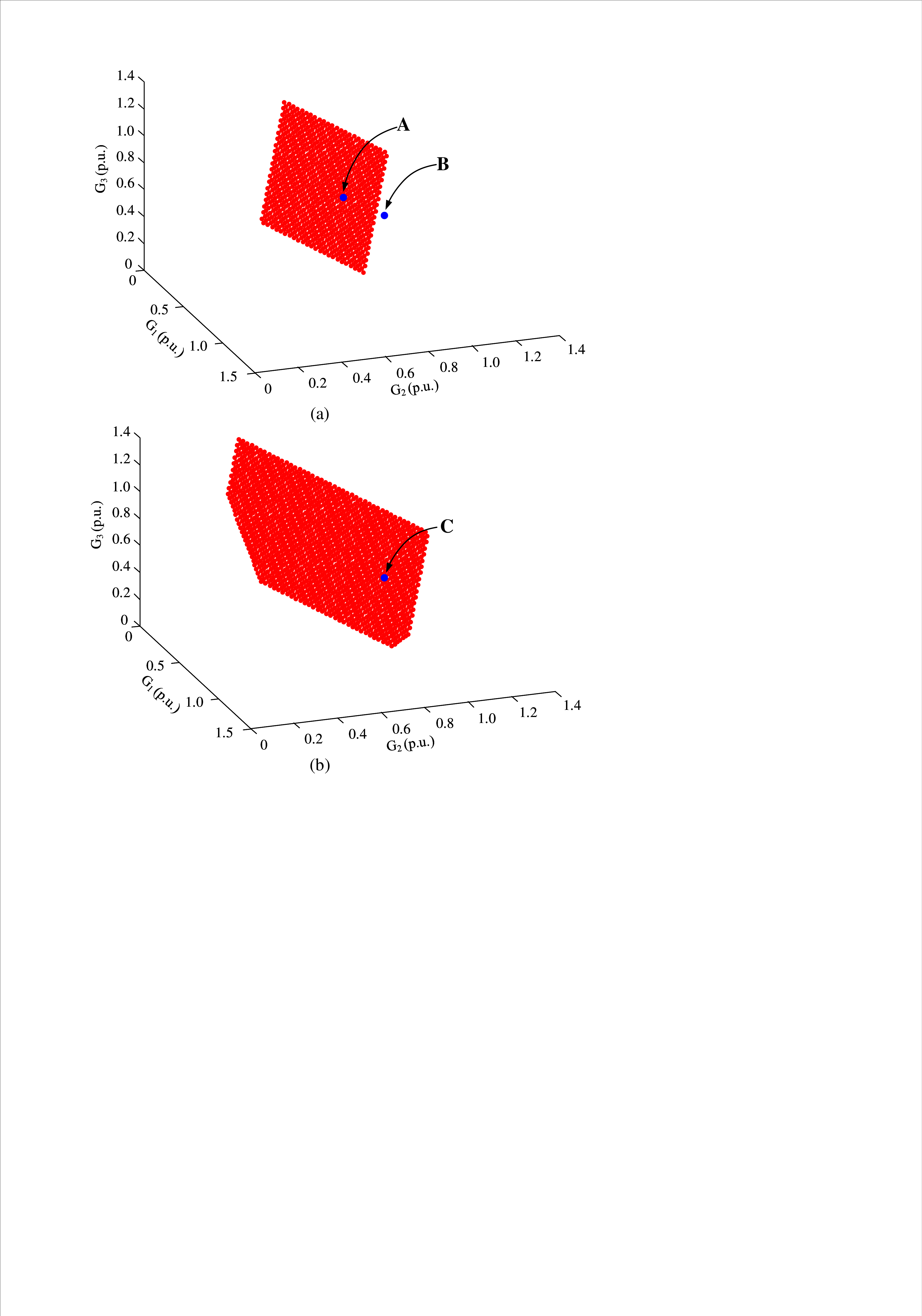}
	\caption{Feasible dispatch region in $t_0$ with ramping rates $0.15$ p.u. (a) and $0.35$ p.u. (b).}
	\label{DifferentRamp_fig}
\end{figure}

\begin{figure}[!h]
	\centering
	\includegraphics[width=7.8cm]{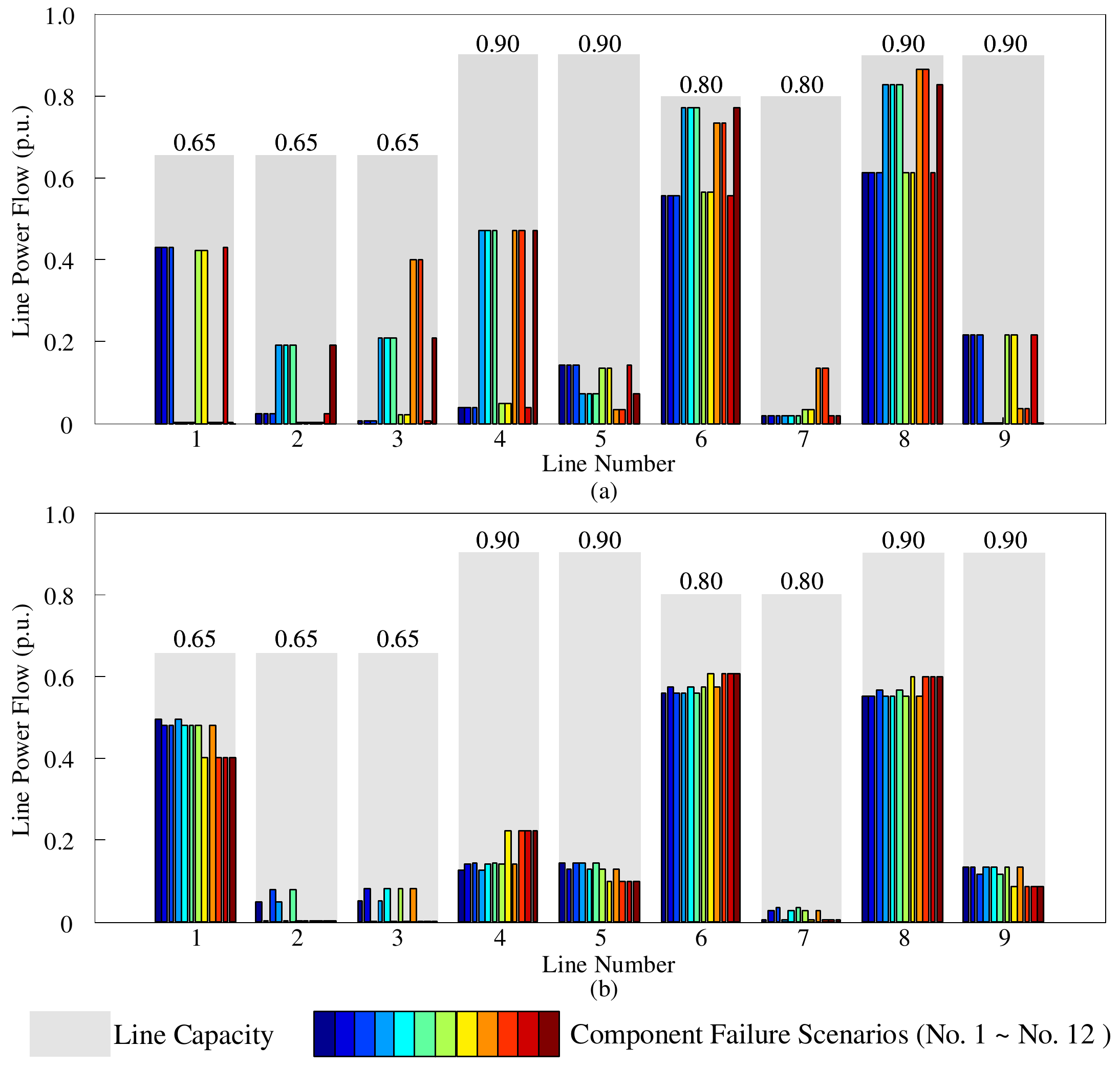}
	\caption{Power flow with ramping rate $0.15$ (p.u.).}
	\label{res15_2540_fig}
\end{figure}

\begin{figure}[!h]
	\centering
	\includegraphics[width=7.8cm]{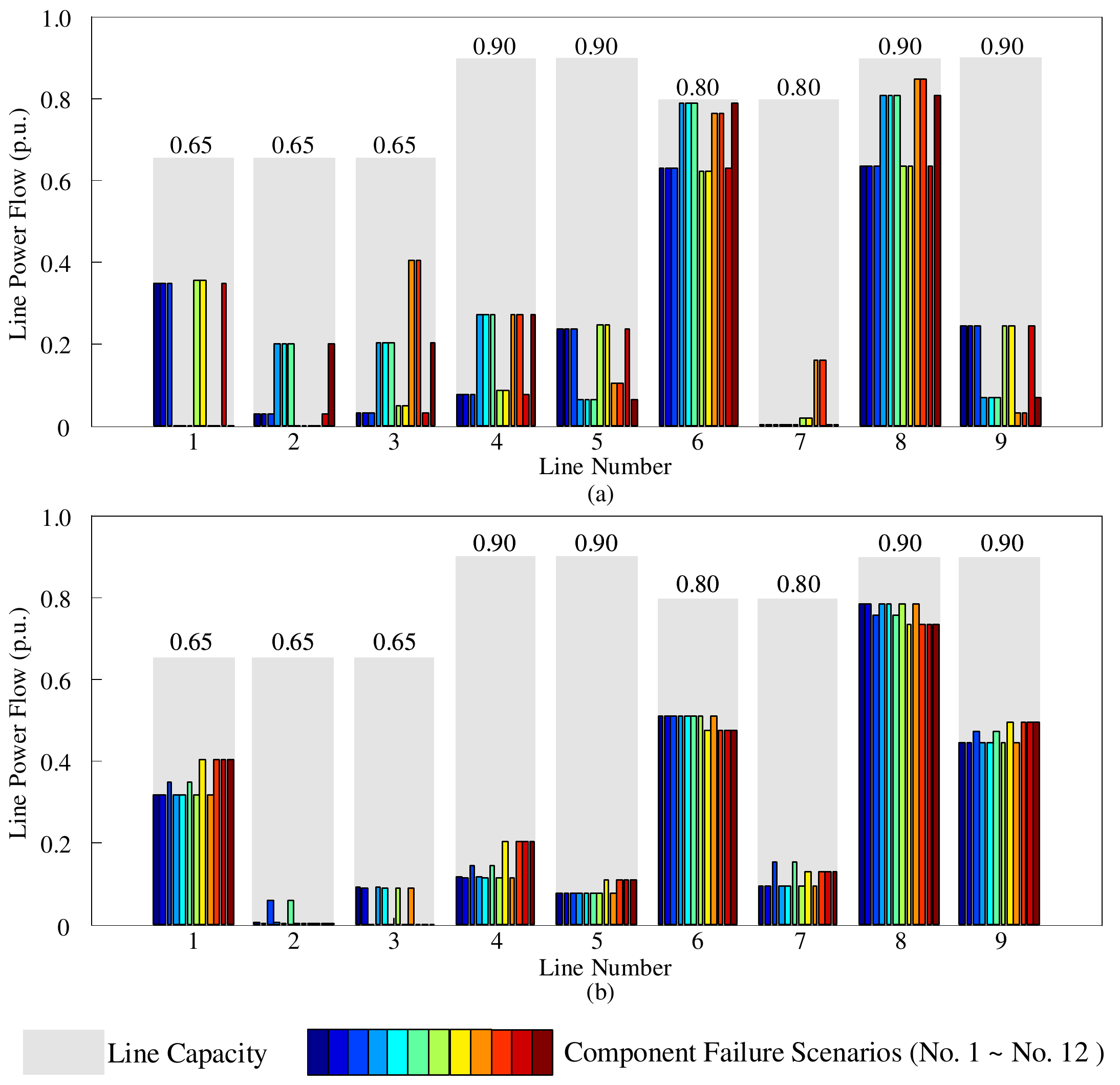}
	\caption{Power Flow with ramping rate $0.35$ (p.u.).}
	\label{res15_2548_fig}
\end{figure}

When making decisions in $t_0$, the possible topology changes in $t_1$ and $t_2$ should be considered, i.e., the extended steady-state security region should be included. Table \ref{IEEE7_PossibleTopo} shows the possible topology changes in $t_1$ and $t_2$. The red regions in Fig. \ref{DifferentRamp_fig} (a) and (b) are the feasible dispatch regions in $t_0$ with the generators' ramping rates $0.15$ (p.u.) and $0.35$ (p.u.), respectively. 
It is observed that the higher generators' ramping rates have a larger feasible dispatch region. The point $A$ corresponds to $P_{G_1}=0.62$, $P_{G_2}=0.70$, and $P_{G_3}=0.68$ in  $t_0$, and this point is within the feasible dispatch region when the generators' ramping rates are $0.15$ (p.u.). Fig. \ref{res15_2540_fig} (a) and (b) show the scenarios of line power flow with possible system topologies caused by the typhnoon in $t_1$ and $t_2$. All scenarios are within the limits. The point $B$ corresponds to $P_{G_1}=0.62$, $P_{G_2}=0.90$, and $P_{G_3}=0.48$ in  $t_0$, and this point is beyond the feasible dispatch region when the generators' ramping rates are $0.15$ (p.u.). Table \ref{IEEE7_PowerFlow} lists line power flow with the component failure scenario $5$. The power flow through the $6^{th}$ line reaches the limit under the optimal strategy in $t_1$ with the loss of load $0.015$ (p.u.) at bus $b_4$. The point $C$ has the same generators' outputs with the point $B$, but it is within the feasible dispatch region when the generators' ramping rates are $0.35$ (p.u.). Fig. \ref{res15_2548_fig} (a) and (b) show the corresponding scenarios of line power flow with possible system topologies caused by the typhnoon in $t_1$ and $t_2$, and all scenarios are within the limits. Higher generators' ramping rates indicates faster adjustment abilities of the system under the operating conditions. 

\begin{figure}[!h]
	\centering
	\includegraphics[width=5.5cm]{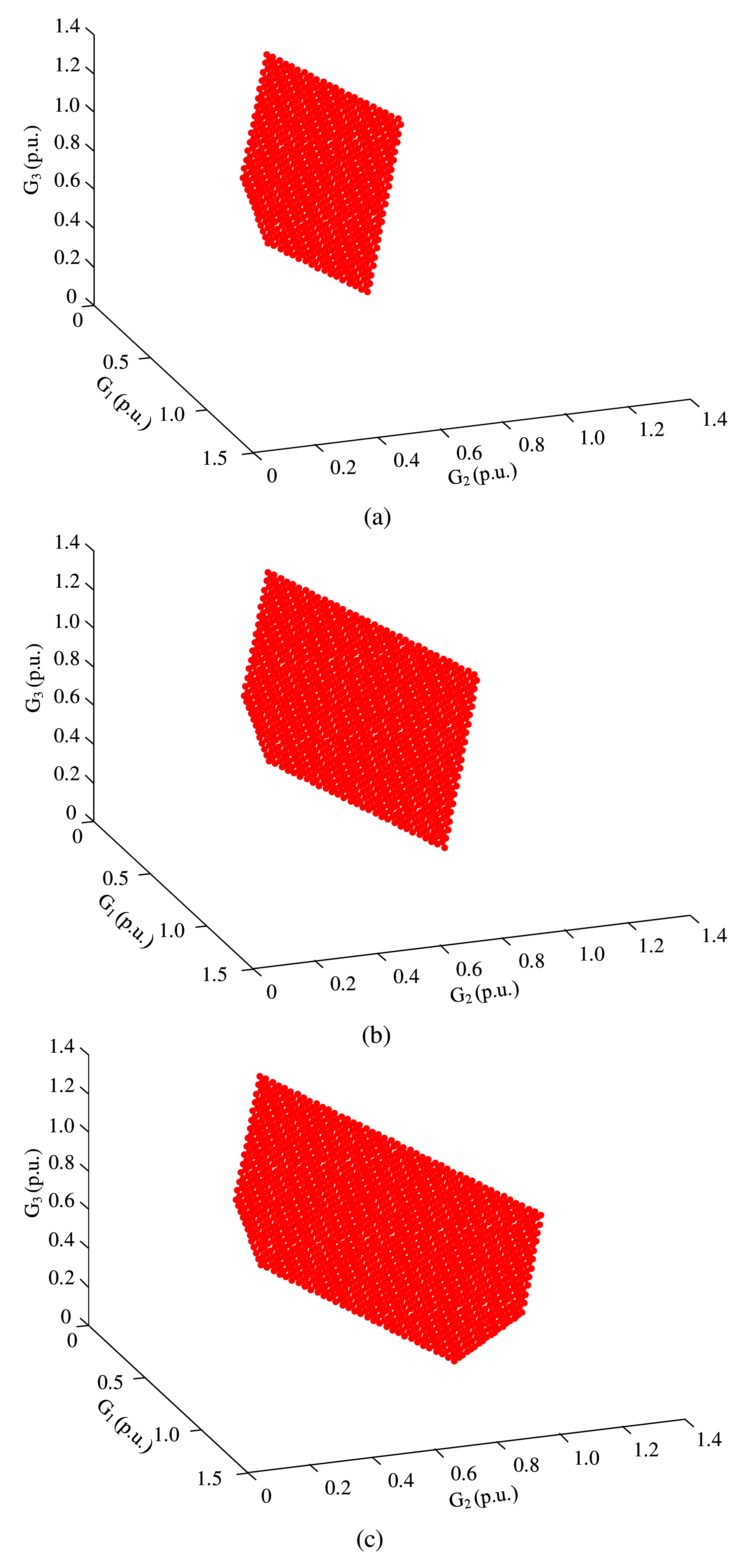}
	\caption{Feasible dispatch regions in $t_0$ with different line capacities $0.7$ (a), $0.8$ (b), $0.9$ (c), respectively.}
	\label{res1_fig}
\end{figure}

Different system parameters also have impacts on the feasible dispatch regions. The generators' ramping rates are set to $0.15$ (p.u.). Fig. \ref{res1_fig} (a), (b), and (c) show the feasible dispatch regions in $t_0$ with different capacities of the line $b_4$-$b_5$, i.e., $0.7$, $0.8$, $0.9$, respectively. It is observed that a larger capacity limit can lead to a larger feasible dispatch region under the same adjustment ability of the system. In practice, the system operators can dispatch the system to the feasible dispatch region to decrease the potential damages. For example, the system operators can dispatch the system to the red region in $t_0$, i.e., the decision period before a typhnoon hits the power system, to enhance system resilience against the typhnoon.   
\begin{table}[!h]
	\renewcommand{\arraystretch}{1.3}
	\caption{Component failure scenarios}
	\label{IEEE7_PowerFlow}
	\begin{tabular}{p{2.2cm}<{\centering}p{2.2cm}<{\centering}p{2.5cm}<{\centering}}
		\hline
		\multirow{2}{*}{No.} & \multicolumn{2}{c}{Line Power Flow}   \\ \cline{2-3} 
		& $t_1$                                                            & $t_2$                  \\ \hline
		1                 & -                                                           & 0.5238    \\ 
		2                 & 0.2050                                                      & -         \\ 
		3                 & 0.1950                                                      & 0.1238    \\ 
		4                 & 0.6050                                                      & 0.2312    \\ 
		5                 & 0.0900                                                      & 0.0345                          \\ 
		6                 & \ct{\fbox{\textbf{0.8000}}}                                                     & 0.5255                          \\ 
		7                 & 0.0100                                                      & 0.0745                           \\ 
		8                 & 0.7850                                                      & 0.5507                          \\ 
		9                 & 0.1150                                                      & 0.1343                          \\   \hline
	\end{tabular}
\end{table}

\subsection{Case 2: IEEE 118-bus System}
\begin{figure}[!h]
	\centering
	\includegraphics[width=8cm]{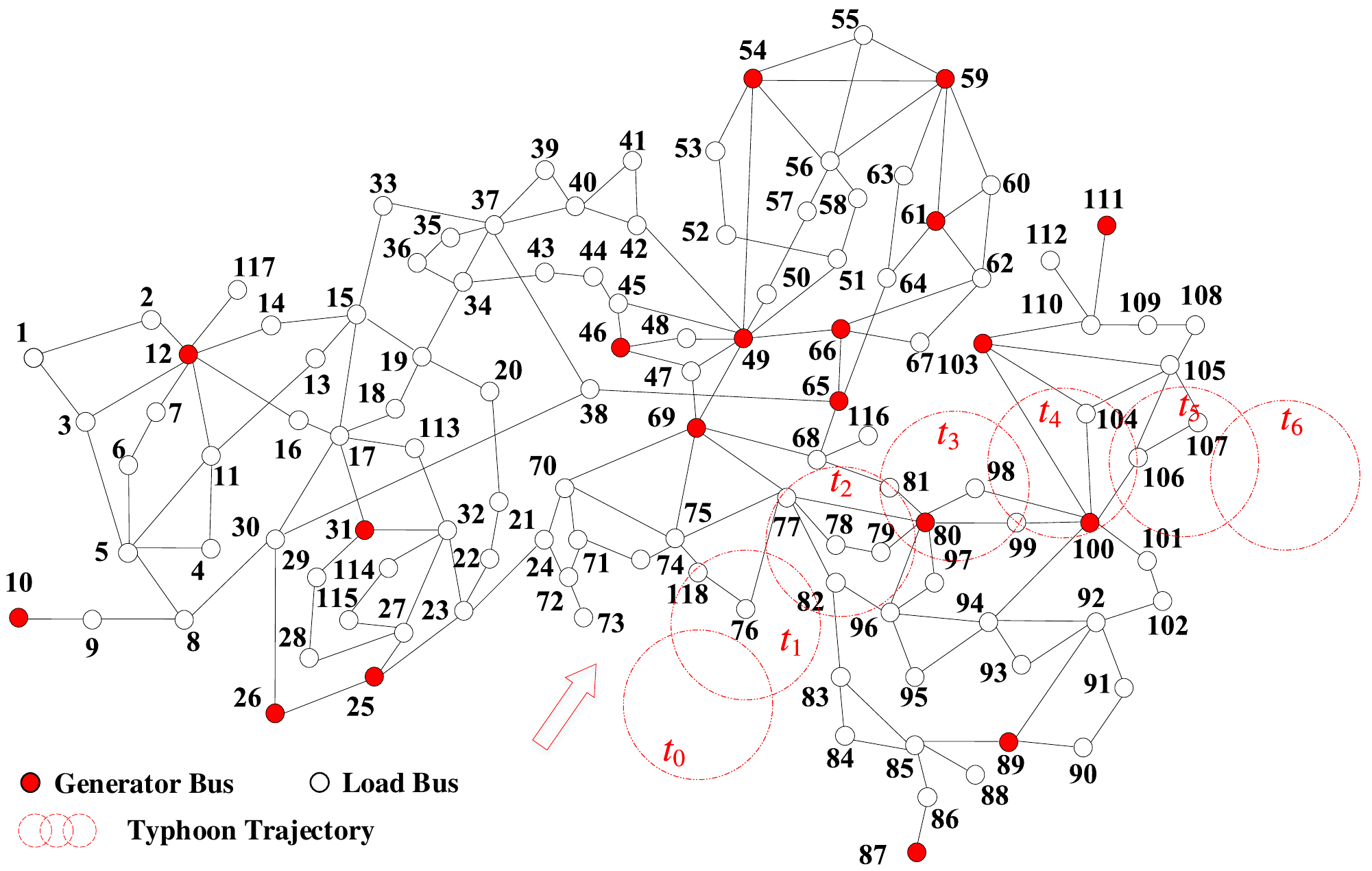}
	\caption{IEEE 118-bus system topology.}
	\label{IEEE118_fig}
\end{figure}

\begin{figure}[!h]
	\centering
	\includegraphics[width=8.5cm]{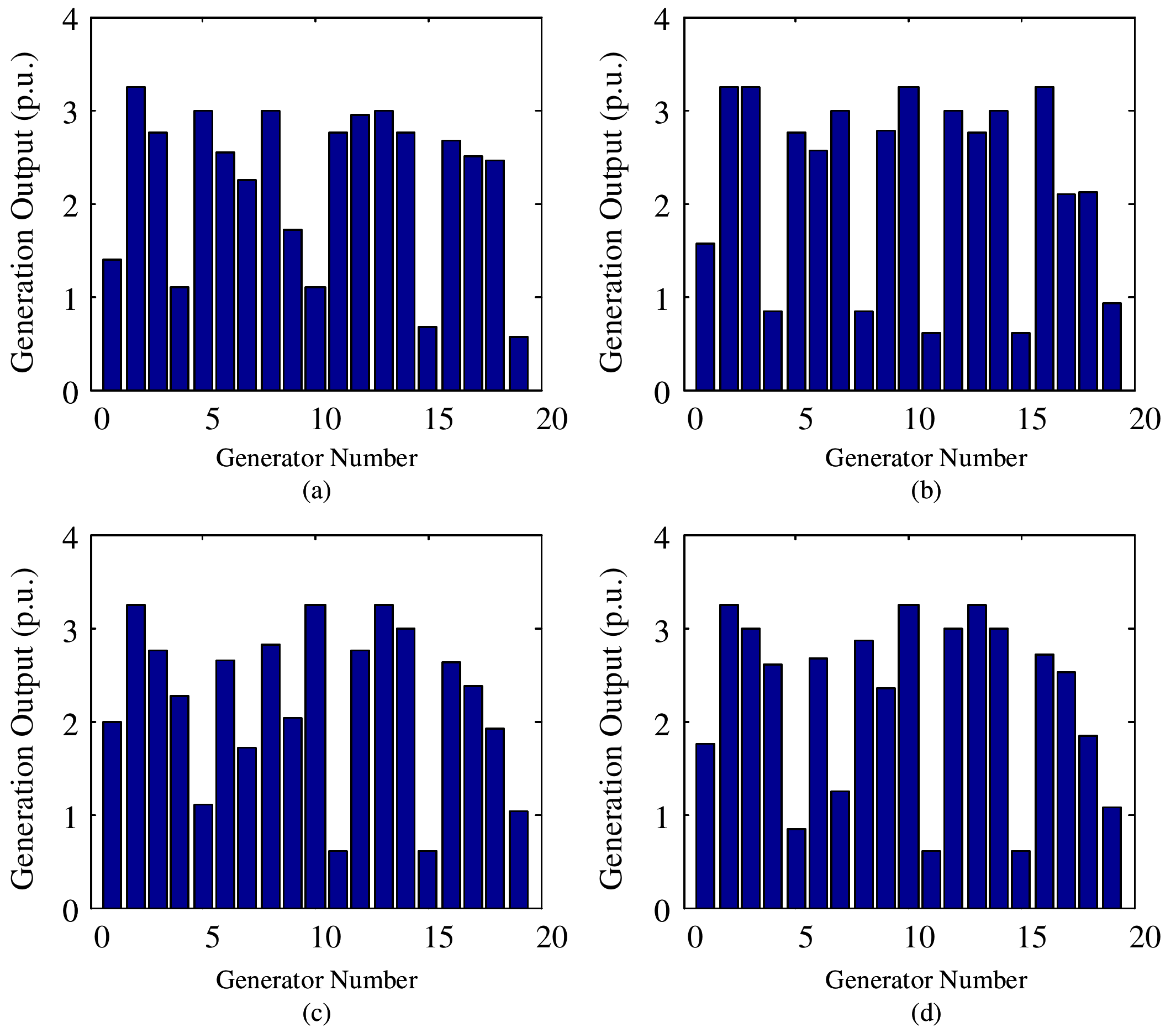}
	\caption{Generation profile with different topology scenarios.}
	\label{IEEE118Res1_fig}
\end{figure}

The IEEE 118-bus System and the trajectory of the typhoon are shown in Fig. \ref{IEEE118_fig}. The ramping rates of all generators are assumed to 0.25. The failure probability of the line on the trajectory of the typhoon is assumed to 0.05. There are numerous uncertain topologies on the typhoon trajectory, and it is difficult to list all uncertain topologies on the typhoon trajectory because of curse of dimensionality with regard to system topologies. In addition, some uncertain topologies on the typhoon trajectory happen with a very small probability. Therefore, the Monte Carlo method is used to generate the topology scenarios. Fig. \ref{IEEE118Res1_fig} (a), (b), (c), and (d) show feasible operating points at time $t_0$ when having $200$, $600$, $1000$, and $1400$ topology scenarios, respectively. It is observed that generation outputs have the similar profiles when considering more topology scenarios such as $1000$ topology scenarios and $1400$ topology scenarios. The similar profiles indicate that the extended steady-state security region based on the topology scenario generation method can approximate the original extended steady-state security region based on all uncertain topology scenarios. For a large system in practice, the system operators first generate enough topology scenarios by using the Monte Carlo method, and then dispatch the system to the feasible dispatch region to against the typhnoon.    

\section{Conclusion}
This paper proposes an extended steady-state security region (ESSR), by which resilient strategies in consideration of uncertain varying topology changes caused by the extreme weather events can be implemented. The physical meaning for the ESSR-based resilient strategy is to find a strategy satisfying the operating constraints in consideration of uncertain varying topology changes. To model this problem, a bilevel optimization model that finds a strategy satisfying the worst-case scenario is established, in which the sub-optimization model in the second level is transformed into the equivalent constraints that are included in the the first level of the model. To address the curse of dimensionality with regard to system topologies for a large-scale system, the Monte Carlo method is used to generate uncertain system topologies, and a McCormick envelope-based approach is proposed to connect generated system topologies to optimization variables. Two test systems are used to validate the proposed model. The detailed results are illustrated in the first test system, in which it is concluded that ESSR-based enhancement strategies are impacted by operating limits such as generators' ramping rates and line thermal capacity limits. The second test system shows the reasonability of the generated system topologies in a large-scale system.

\bibliographystyle{IEEEtran}
\bibliography{IEEEabrv,RefDatabase}


\end{document}